\documentclass[12pt]{article}
\usepackage{amssymb}

\newtheorem{definition}{Definition}[section]
\newtheorem{theorem}[definition]{Theorem}
\newtheorem{lemma}[definition]{Lemma}
\newtheorem{corollary}[definition]{Corollary}
\newtheorem{remark}[definition]{Remark}

\newtheorem{note}[definition]{Note}
\newtheorem{notation}[definition]{Notation}
\newtheorem{proposition}[definition]{Proposition}

\def\qed{{~~~\vrule height .75em width .4em depth .2em}}

\def\Z{\mathbb Z}
\def\I{\mathbb I}
\def\O{\mathcal{O}}

\newcommand{\beast}{\begin{eqnarray*}}
\newcommand{\eeast}{\end{eqnarray*}}

\begin{document}

\title{ \bf The Tetrahedron Algebra and its Finite-Dimensional
  Irreducible Modules}

\author{Brian Hartwig {\footnote{
Department of Mathematics, University of
Wisconsin, 480 Lincoln Drive, Madison WI 53706-1388 USA}
}}
\date{}
%to get date printout, comment out above line
\maketitle

\begin{abstract}
Recently Terwilliger and the present author found a
presentation for the three-point $\mathfrak{sl}_2$ loop algebra via
generators and relations.  To obtain this presentation we defined
a Lie algebra $\boxtimes$ by generators and relations and
displayed an isomorphism from $\boxtimes$ to the three-point
$\mathfrak{sl}_2$ loop algebra.  In this paper we classify the
finite-dimensional irreducible $\boxtimes$-modules.
\end{abstract}

\section{Introduction and Statement of Results}

In \cite{HandT} Terwilliger and the present author gave a
presentation of the three-point $\mathfrak{sl}_2$ loop algebra via
generators and relations.  To obtain this presentation we defined
a Lie algebra $\boxtimes$ by generators and relations and
displayed an isomorphism from $\boxtimes$ to the three-point
$\mathfrak{sl}_2$ loop algebra.

In this paper we classify the finite-dimensional
irreducible $\boxtimes$-modules.  To obtain this classification, we
exploit a connection between $\boxtimes$ and a certain
infinite-dimensional Lie algebra which appears in the physics literature,
called the Onsager algebra \cite{A}, \cite {DG}, \cite{GR}, \cite{ons1}.
Before we explain our results we first give some
background on $\boxtimes$ and the Onsager algebra.  We start with their
definitions.

Throughout the paper $\Bbb K$ denotes an algebraically closed field
of characteristic $0$.

\begin{definition}
\label{def:main} \rm \cite{HandT}
Let $\boxtimes$ denote the Lie algebra over $\Bbb K$ with generators
$$ \{ X_{rs} | r,s \in \I, r \not=s \} \qquad \qquad \I = \{0,1,2,3 \}
$$
and the following relations.
\begin{enumerate}
\item
For all distinct $r,s \in \I$,
\begin{eqnarray}
X_{rs}+ X_{sr} =0. \label{eq:1}
\end{eqnarray}
\item
For all mutually distinct $r,s,t \in \I$,
\begin{eqnarray}
[X_{rs},X_{st}]=2X_{rs} + 2X_{st}. \label{eq:2}
\end{eqnarray}
\item
For all mutually distinct $r,s,t,u \in \I$,
\begin{eqnarray}
[X_{rs},[X_{rs},[X_{rs},X_{tu}]]] = 4[X_{rs},X_{tu}]. \label{eq:3}
\end{eqnarray}
\end{enumerate}
We call $\boxtimes$ the {\it tetrahedron algebra}.
\end{definition}

\begin{definition}
\rm \cite{ons1}, \cite{P} \label{def:oalg} Let $\O$ denote the Lie
algebra over $\Bbb K$ with generators $X,Y$ and relations
\begin{eqnarray}
[X,[X,[X,Y]]] = 4 [X,Y], \label{eq:fo1} \\
\left[ Y,[Y,[Y,X]] \right] = 4 [Y,X]. \label{eq:fo2}
\end{eqnarray}
We call $\cal O$ the {\it Onsager algebra}. We call $X,Y$ the {\it
standard generators} for $\O$.
\end{definition}

We recall the connection between $\boxtimes$ and $\O$.  For mutually
distinct $r,s,t,u \in \I$ there exists a Lie algebra
injection from $\O$ into $\boxtimes$ that sends \beast
X & \mapsto & X_{rs} \\
Y & \mapsto & X_{tu} \eeast \cite[Corollary 12.2]{HandT}.  We
call the image of this injection an {\it Onsager subalgebra} of
$\boxtimes$.  Observe that $\boxtimes$ has three Onsager subalgebras.
By \cite[Theorem 11.6]{HandT} the $\Bbb K$-vector space $\boxtimes$ is the
direct sum of its three Onsager subalgebras.

We now summarize the classification of finite-dimensional irreducible
$\O$-modules \cite{Dav1}, \cite{Dav2}.  We begin
with a comment.  Let $V$ denote a finite-dimensional irreducible
$\O$-module.  As we will see in Section \ref{sect:2}, the standard
generators $X,Y$ are diagonalizable on $V$.  Furthermore there
exist an integer $d \geq 0$ and scalars $\alpha, \alpha ^* \in
\Bbb K$ such that the set of distinct eigenvalues of $X$ (resp.
$Y$) on $V$ is $\{d-2i+\alpha | 0 \leq i \leq d \}$ (resp.
$\{d-2i+\alpha^* | 0 \leq i \leq d \}$).  We call the ordered pair
$(\alpha , \alpha^*)$ the {\it type} of $V$.  Replacing $X,Y$ by
$X-\alpha I, Y-\alpha ^* I$ the type becomes $(0,0)$.  Therefore it
suffices to classify the finite-dimensional irreducible $\O$-modules
of type $(0,0)$.

We begin with a special case.  Observe that up to isomorphism, there
exists a unique irreducible $\O$-module of dimension $1$ and type
$(0,0)$.  We call this the {\it trivial} $\O$-module.

Let $\mathfrak{sl}_2$ denote the Lie algebra over $\Bbb K$ with basis
$e,f,h$ and Lie bracket $$ \left[ e,f \right] = h, \qquad \left[ h,e
\right] = 2e, \qquad
\left[ h,f \right] = -2f. $$ Let $t$ denote an indeterminate and let
$L(\mathfrak{sl}_2)$ denote the loop algebra $\mathfrak{sl}_2 \otimes
\Bbb K [t,t^{-1}]$, where $\otimes$ means $\otimes_{\Bbb K}$.  By
\cite[Proposition 1]{roan} there exists a Lie algebra injection from
$\O$ into $L(\mathfrak{sl}_2)$ that sends \beast
X & \mapsto & e \otimes 1 + f \otimes 1 \\
Y & \mapsto & e \otimes t + f \otimes t^{-1}. \eeast For the moment
we identify $\O$ with its image under the above injection.  For
nonzero $a \in \Bbb K$ we define the Lie algebra homomorphism
$EV_a: L(\mathfrak{sl}_2) \rightarrow \mathfrak{sl}_2$ by $EV_a(u
\otimes g(t)) = g(a)u$ for all $u \in \mathfrak{sl}_2$ and $g(t) \in \Bbb
K [t,t^{-1}]$.  Let $ev_a$ denote the restriction of
$EV_a$ to $\O$.  Then $ev_a:\O \rightarrow \mathfrak{sl}_2$ is a Lie
algebra homomorphism which we call the {\it evaluation homomorphism}
for $a$. Let $V$ denote an irreducible $\mathfrak{sl}_2$-module with
finite dimension at least $2$.  We pull back by $ev_a$ to get
an $\O$-module structure on $V$.  We call this an {\it evaluation
module} for $\O$ and denote it by $V(a)$.  The $\O$-module $V(a)$
is irreducible if and only if $a \not= \pm 1$ \cite[Lemma 4]{DR}.
Irreducible $\O$-modules $V(a)$ and $V(b)$ are isomorphic if and only
if $a = b^{\pm 1}$ \cite[Proposition 5]{roan}.

Let $U,V$ denote $\O$-modules.  Then $U \otimes V$
has an $\O$-module structure given by
$$ x.(u \otimes v) = (x.u) \otimes v + u \otimes (x.v) \qquad \quad
x \in \O, \quad u \in U, \quad v \in V.$$ Let $V$ denote an
$\O$-module that is the tensor product of finitely many evaluation
modules.  If $V$ is irreducible then it is type $(0,0)$ \cite[p.
3281]{DR}.

The classification of finite-dimensional irreducible $\O$-modules
of type $(0,0)$ is given in the following three theorems.

\begin{theorem} \cite[Theorem 6]{DR}
Every nontrivial finite-dimensional irreducible $\O$-module of type
$(0,0)$ is isomorphic to a tensor product of evaluation modules.
\end{theorem}

\begin{theorem} \cite[Proposition 5]{DR}
Let $V_1(a_1), \ldots, V_n(a_n)$ denote a finite sequence of
evaluation modules for $\O$, and consider the $\O$-module $V_{1}(a_1)
\otimes \cdots \otimes V_{n}(a_n)$.  This module is irreducible if and
only if $a_1, a_1^{-1}, \ldots, a_n,a_n^{-1}$ are
mutually distinct.
\end{theorem}

\begin{definition} \rm
Let $V_1(a_1), \ldots, V_n(a_n)$ denote a finite sequence of
evaluation modules for $\O$.  Let $V$ denote the $\O$-module
$V_{1}(a_1) \otimes \cdots \otimes V_{n}(a_n)$.  Any tensor product of
evaluation modules that can be obtained from $V$ by permuting the
order of the factors and replacing any number of the $a_i$'s with
their multiplicative inverses will be called {\it equivalent} to $V$.
\end{definition}

\begin{theorem}  \cite[Proposition 5]{DR}
Let $U$ and $V$ denote tensor products of finitely many evaluation modules
for $\O$.  Assume each of $U,V$ is irreducible as an $\O$-module.
Then the $\O$-modules $U$ and $V$ are isomorphic if and only if they
are equivalent.
\end{theorem}

This completes the classification of finite-dimensional irreducible
$\O$-modules. 

We now state our main results.  They are contained in the
following two theorems and subsequent remark.

\begin{theorem} \label{th:H1}  Let $V$ denote a finite-dimensional
irreducible $\boxtimes$-module.  Then there exists a unique
$\O$-module structure on $V$ such that the standard generators $X,Y$ act
on $V$ as $X_{01},X_{23}$ respectively.  This $\O$-module structure is
irreducible and has type $(0,0)$.
\end{theorem}

\begin{theorem}  \label{th:H2}  Let $V$ denote a
finite-dimensional irreducible $\O$-module of type $(0,0)$. Then
there exists a unique $\boxtimes$-module structure on $V$ such that
the standard generators $X,Y$ act on $V$ as $X_{01},X_{23}$
respectively. This $\boxtimes$-module structure is irreducible.
\end{theorem}

\begin{remark}
Combining the previous two theorems we obtain a bijection between
the following two sets:
\begin{enumerate}
\item
The isomorphism classes of finite-dimensional irreducible
$\O$-modules of type $(0,0)$.
\item
The isomorphism classes of finite-dimensional irreducible
$\boxtimes$-modules.
\end{enumerate}
\end{remark}

\section{$\O$-Modules and Tridiagonal Pairs \label{sect:2}}

In order to prove Theorem \ref{th:H1} and Theorem \ref{th:H2} it will
be useful to consider how finite-dimensional irreducible $\O$-modules
are related to tridiagonal pairs.  To explain this relationship we
use the following concepts.  Let $V$ denote a vector space
over $\Bbb K$ with positive finite dimension.  By a {\it linear
transformation on $V$} we mean a $\Bbb K$-linear map from $V$ to
$V$.  Let $A$ denote a linear transformation on $V$.  For $\lambda
\in \Bbb K$ define
\begin{eqnarray}
V_A(\lambda) := \{ v \in V | A v = \lambda v \}.  \label{eq:espace}
\end{eqnarray}
Observe that $\lambda$ is an eigenvalue for $A$ if and only if
$V_A(\lambda) \not= 0$, and in this case $V_A(\lambda)$ is the
corresponding eigenspace of $A$.
We recall that the sum $\sum_{\lambda \in \Bbb
K}V_A(\lambda)$ is direct.  We call $A$ {\it diagonalizable} whenever
$V = \sum_{\lambda \in \Bbb K}V_A(\lambda)$.

\begin{lemma}
\label{lem:conf}
Let $V$ denote a vector space over $\Bbb K$ with positive finite
dimension.  Let $A$ and $A^*$ denote linear transformations on $V$.
Then for $\lambda \in \Bbb K$ the following $(i),(ii)$ are equivalent.
\begin{enumerate}
\item
The expression $[A,[A,[A,A^*]]] - 4[A,A^*]$ vanishes on
$V_A(\lambda)$.
\item
\begin{eqnarray}
\label{eq:un2}
A^* V_A( \lambda) \subseteq V_A( \lambda+2)+V_A( \lambda)+V_A( \lambda-2).
\end{eqnarray}
\end{enumerate}
\end{lemma}

{\it Proof.}
Let $\Phi$ denote the expression in $(i)$ and observe
\beast
\Phi = A^3 A^* - 3A^2 A^* A + 3A A^* A^2 - A^* A^3 - 4A
A^* + 4 A^* A.
\eeast
For $v \in V_A(\lambda)$ we evaluate $\Phi v$ using $Av = \lambda v$
to find
\beast
\Phi v &=& (A^3 A^* - 3 \lambda A^2 A^* + 3 \lambda^2 A A^* -
\lambda^3 A^* - 4A A^* + 4 \lambda A^*)v \\
&=& (A - (\lambda+2)I)(A - \lambda I)(A - (\lambda-2)I)A^*v.
\eeast The scalars $\lambda+2,\lambda,\lambda-2$ are mutually
distinct since ${\rm Char} (\Bbb K) = 0$.  The result follows.
\qed

\medskip

We now recall the concept of a {\it tridiagonal pair}
\cite{bandt}, \cite{iandt}, \cite{ITT}, \cite{ito}.

\begin{definition}
\label{def:tp} \rm \cite{ITT} Let $V$ denote a vector space over
$\Bbb K$ with positive finite dimension.  By a {\it tridiagonal
pair} on $V$, we mean an ordered pair $A,A^*$ where $A$ and $A^*$
are linear transformations on $V$ that satisfy the following four
conditions.
\begin{enumerate}
\item Each of $A,A^*$ is diagonalizable.
\item There exists an ordering $V_0,V_1, \ldots, V_d$ of the
eigenspaces of $A$ such that
\begin{eqnarray}
\label{eq:un1}
A^* V_i \subseteq V_{i+1} + V_i + V_{i-1} \qquad \qquad 0 \leq i \leq d,
\end{eqnarray}
where $V_{-1}=0$, $V_{d+1} =0$.
\item  There exists an ordering $V_0^*,V_1^*, \ldots, V_{\delta}^*$ of the
eigenspaces of $A^*$ such that
$$A V_i^* \subseteq V_{i+1}^* + V_i^* + V_{i-1}^* \qquad \qquad 0 \leq i
\leq \delta,$$
where $V_{-1}^*=0$, $V_{\delta+1}^* =0$.
\item  There does not exist a subspace $W$ of $V$ such that $AW
\subseteq W$, $A^*W \subseteq W$, $W \not= 0$, and $W \not= V$.
\end{enumerate}
\end{definition}

We comment on Definition \ref{def:tp}.  Let $A,A^*$ denote a
tridiagonal pair on $V$ and let $d$ and $\delta$ be as in
Definition \ref{def:tp}.  From \cite[Lemma 4.5]{ITT} we find $d=
\delta$.  We call this common value the {\it diameter} of $A,A^*$.
An ordering of the eigenspaces of $A$ (resp. $A^*$) will be called
{\it standard} whenever it satisfies Definition \ref{def:tp}(ii)
(resp. Definition \ref{def:tp}(iii)).   Let $V_0, \ldots, V_d$
denote a standard ordering of the eigenspaces of $A$.  By \cite[p.
388]{ter2}, the ordering $V_d, \ldots, V_0$ is standard
and no other ordering is standard.  A similar result holds for
$A^*$.  By an {\it eigenvalue sequence} (resp. {\it dual
eigenvalue sequence}) of $A,A^*$ we mean an ordering of the
eigenvalues for $A$ (resp. $A^*$) for which the corresponding
ordering of the eigenspaces is standard.

\begin{definition}
\rm Let $\theta_0, \theta_1, \ldots, \theta_d$ denote a finite
sequence of scalars in $\Bbb K$.  We say this sequence is {\it
arithmetic with common difference $2$} whenever $\theta_{i-1} -
\theta_{i}= 2$ for $1 \leq i \leq d$.
\end{definition}

\begin{theorem}
\label{th:p1} Let $V$ denote a vector space over $\Bbb K$ with
positive finite dimension.  Let $A$ and $A^*$ denote linear
transformations on $V$.  Then the following $(i),(ii)$ are
equivalent.
\begin{enumerate}
\item
$A,A^*$ is a tridiagonal pair on $V$ whose eiqenvalue
sequence and dual eigenvalue sequence are both
arithmetic with common difference $2$.
\item
There exists an irreducible $\O$-module structure on $V$ such that
the standard generators $X,Y$ act on $V$ as $A,A^*$ respectively.
\end{enumerate}
\end{theorem}

{\it Proof.} $(i) \Rightarrow (ii):$  We
first show that $A,A^*$ satisfy
\begin{eqnarray}
\label{eq:o1} [A,[A,[A,A^*]]] = 4[A,A^*].
\end{eqnarray}
By assumption the tridiagonal pair $A,A^*$ has an eigenvalue
sequence that is arithmetic with common difference $2$. Let
$V_0, V_1, \ldots, V_d$ denote the corresponding ordering of the
eigenspaces of $A$.  This ordering is standard by construction so
it satisfies (\ref{eq:un1}).  By these comments we find that $A^*$
satisfies (\ref{eq:un2}) for all $\lambda \in \Bbb K$.  By Lemma
\ref{lem:conf} and since $A$ is diagonalizable on $V$ we find $A,A^*$
satisfy (\ref{eq:o1}).  Reversing the roles of $A$ and $A^*$ in
the above argument we find $A,A^*$ satisfy
\begin{eqnarray}
\label{eq:o2} [A^*,[A^*,[A^*,A]]] = 4[A^*,A].
\end{eqnarray}
By (\ref{eq:o1}) and (\ref{eq:o2}) there exists an $\O$-module
structure on $V$ such that $X,Y$ act as $A,A^*$ respectively.
This $\O$-modules structure is irreducible by Definition
\ref{def:tp}(iv).

$(ii) \Rightarrow (i):$  \cite[Example 1.6]{ITT} Since $\Bbb K$ is
algebraically closed
all the eigenvalues for $A$ are contained in $\Bbb K$.  Since $V$
has positive finite dimension, $A$ has at least one eigenvalue
$\lambda$.  Since ${\rm Char} (\Bbb K) = 0$, the scalars $\lambda,
\lambda+2, \ldots$ are mutually distinct and therefore
cannot all be eigenvalues for $A$.  Hence there exists an
eigenvalue $\theta$ for $A$ such that $\theta +2$ is not an
eigenvalue for $A$.  The scalars $\theta, \theta-2,
\ldots$ are mutually distinct and therefore cannot all be
eigenvalues for $A$.  Hence there exists a nonnegative integer $d$
such that $\theta - 2i$ is an eigenvalue of $A$ for $0 \leq i
\leq d$ but is not an eigenvalue of $A$ for $i = d+1$. Abbreviate
$V_A(\theta -2i)$ by $V_i$ for $0 \leq i \leq d$.  By construction
$\sum_{i=0}^{d} V_i$ is $A$-invariant. By Lemma \ref{lem:conf} we
find $$A^* V_i \subseteq V_{i+1} + V_i + V_{i-1} \qquad (0 \leq i
\leq d),$$ where $V_{-1}=0$, $V_{d+1} =0$. Therefore
$\sum_{i=0}^{d} V_i$ is $A^*$-invariant.  Since the $\O$-module
structure is irreducible and since $\sum_{i=0}^{d} V_i \not= 0$ we
find $\sum_{i=0}^{d} V_i = V$.  We have now shown that $A$ is
diagonalizable and Definition \ref{def:tp}(ii) holds. Reversing
the roles of $A$ and $A^*$ in the above argument we find $A^*$ is
diagonalizable and Definition \ref{def:tp}(iii) holds. Definition
\ref{def:tp}(iv) is immediate since the $\O$-module $V$ is
irreducible. We have now shown that $A,A^*$ is a tridiagonal pair
on $V$. Recall that $\theta -2i$ is the eigenvalue for $A$
associated with $V_i$ for $0 \leq i \leq d$. The ordering
$V_0, \ldots, V_d$ is standard so the sequence $\theta -2i$
$(0 \leq i \leq d)$ is an eigenvalue sequence for $A,A^*$. This
sequence is arithmetic with common difference $2$.  We have now
shown $A,A^*$ has an arithmetic eigenvalue sequence with common
difference $2$. Reversing the roles of $A$ and $A^*$ in the above
argument we find that $A,A^*$ has an arithmetic dual eigenvalue
sequence with common difference $2$. \qed

\medskip

\begin{definition}
\rm Let $V$ denote a finite-dimensional irreducible $\O$-module.
By Theorem \ref{th:p1} there exist an integer $d \geq 0$ and
scalars $\alpha, \alpha ^* \in \Bbb K$ such that the set of
distinct eigenvalues of $X$ (resp. $Y$) on $V$ is $\{d-2i+\alpha |
0 \leq i \leq d \}$ (resp. $\{d-2i+\alpha^* | 0 \leq i \leq d
\}$). We call $d$ the {\it diameter} of $V$.  We call the ordered pair
$(\alpha, \alpha^*)$ the {\it
type} of $V$.
\end{definition}

\begin{note}
\rm Let $V$ denote a finite-dimensional irreducible $\O$-module of
type $(\alpha, \alpha^*)$.  Replacing $X,Y$ by $X-\alpha I,
Y-\alpha ^* I$ the type becomes $(0,0)$.
\end{note}

Restricting Theorem \ref{th:p1} to type $(0,0)$ we get the
following corollary.

\begin{corollary}
\label{cor:t0o} Let $V$ denote a vector space over $\Bbb K$ with
positive finite dimension.  Let $A$ and $A^*$ denote linear
transformations on $V$.  Then the following $(i),(ii)$ are
equivalent.
\begin{enumerate}
\item
$A,A^*$ is a tridiagonal pair on $V$ and $d-2i$ $(0 \leq i \leq
d)$ is both an eigenvalue sequence and dual eigenvalue sequence
for $A,A^*$, where $d$ denotes the diameter.
\item
There exists an irreducible $\O$-module structure on $V$ of type
$(0,0)$ such that the standard generators $X,Y$ act on $V$ as $A,A^*$
respectively.
\end{enumerate}
\end{corollary}

{\it Proof.}  Immediate from Theorem \ref{th:p1} and the definition of
type.

\medskip

\section{Finite-Dimensional Irreducible $\boxtimes$-Modules}

Let $V$ denote a finite-dimensional irreducible
$\boxtimes$-module. In this section we obtain the following
description of $V$. We show that each generator $X_{rs}$ of
$\boxtimes$ is diagonalizable on $V$.  We show that there exists
an integer $d \geq 0$ such that for each $X_{rs}$ the eigenvalues
for $X_{rs}$ on $V$ are $d, d-2, \ldots, -d$.  We give the action
of each $X_{rs}$ on the eigenspaces of the other generators.  In
our investigation we use the following lemma.

\begin{lemma}
\label{lem:conf2}
Let $V$ denote a vector space over $\Bbb K$ with positive finite
dimension.  Let $A$ and $B$ denote linear transformations on $V$.
Then for $\lambda \in \Bbb K$ the following $(i),(ii)$ are
equivalent.
\begin{enumerate}
\item
The expression $[A,B] - 2A - 2B$ vanishes on $V_A(\lambda)$.
\item
$(B + \lambda I) V_A(\lambda) \subseteq V_A(\lambda+2).$
\end{enumerate}
\end{lemma}

{\it Proof.}
Let $\Psi$ denote the expression in $(i)$ and observe
\beast
\Psi = AB-BA-2A-2B.
\eeast
For $v \in V_A(\lambda)$ we evaluate $\Psi v$ using $Av= \lambda v$
and find
\beast
\Psi v &=& (A - (\lambda + 2)I)(B + \lambda I)v.
\eeast
The result follows.
\qed

\medskip

We refine notation (\ref{eq:espace}) as follows.

\begin{notation}
\rm \label{note:u} With reference to Definition \ref{def:main} let
$V$ denote a finite-dimensional irreducible $\boxtimes$-module.
For distinct $r,s \in \I$ and for all $\lambda \in \Bbb K$ we
define
\begin{eqnarray} \label{eq:es2} V_{rs} (\lambda) = \{ v \in V |
X_{rs} v = \lambda v \}.
\end{eqnarray}
\end{notation}

\begin{theorem}
\label{lem:1} Let $V$ denote a finite-dimensional irreducible
$\boxtimes$-module.  For all $r,s,t,u \in \I$ ($r \not= s, t \not=
u$) and for all $\lambda \in \Bbb K$ the action of $X_{tu}$ on
$V_{rs}(\lambda)$ is given as follows.

\medskip

\begin{tabular}{c|c}
Case & Action of $X_{tu}$ on $V_{rs}(\lambda)$ \\
\hline
$t=r, \qquad u=s$ & $(X_{tu} - \lambda I) V_{rs}(\lambda ) = 0$ \\
$t=s, \qquad u=r$ & $(X_{tu} + \lambda I) V_{rs}(\lambda ) = 0$ \\
$t=s, \qquad u \not= r$ & $(X_{tu} + \lambda I) V_{rs}(\lambda ) \subseteq V_{rs}(\lambda + 2)$ \\
$t \not= r, \qquad u=s$ & $(X_{tu} - \lambda I) V_{rs}(\lambda )
\subseteq V_{rs}(\lambda + 2)$ \\
$t=r, \qquad u \not= s$ & $(X_{tu} - \lambda I) V_{rs}(\lambda )
\subseteq V_{rs}(\lambda - 2)$ \\
$t \not= s, \qquad u=r$ & $(X_{tu} + \lambda I) V_{rs}(\lambda )
\subseteq V_{rs}(\lambda - 2)$ \\
$r,s,t,u$ distinct & $X_{tu} V_{rs}(\lambda ) \subseteq
V_{rs}(\lambda + 2) + V_{rs}(\lambda ) + V_{rs}(\lambda - 2) $
\end{tabular}

\medskip

\noindent
 We are using Notation \ref{note:u}.
\end{theorem}

{\it Proof.}  We consider each row in the above table.

$t=r, u=s$: Immediate from (\ref{eq:es2}).

$t=s, u=r$: Immediate from (\ref{eq:es2}) and since $X_{rs}
= - X_{tu}$ by Definition \ref{def:main}(i).

$t=s, u \not= r$: Combine Definition \ref{def:main}(ii) and
Lemma \ref{lem:conf2}.

$t \not= r, u=s$: By case $t=s, u \not= r$ above and since
$X_{ut} = -X_{tu}$.

$t=r, u \not= s$: By case $t=s,u \not= r$ above and since
$V_{rs}(\lambda) = V_{sr}(-\lambda)$.

$t \not= s, u=r$: By case $t=r, u \not= s$ above and since $X_{ut}
= -X_{tu}$.

$r,s,t,u$ distinct: Combine Definition \ref{def:main}(iii) and
Lemma \ref{lem:conf}. \qed

\begin{corollary}
\label{cor:act} Let $V$ denote a finite-dimensional irreducible
$\boxtimes$-module.  Then for all $r,s,t,u \in \I$ ($r \not= s, t
\not= u$) and for all $\lambda \in \Bbb K$ we have
$$X_{tu} V_{rs}(\lambda) \subseteq V_{rs}(\lambda+2) + V_{rs}(\lambda) +
V_{rs}(\lambda-2).$$ We are using Notation \ref{note:u}.
\end{corollary}

{\it Proof.}
The above inclusion holds for each row of the table in Theorem
\ref{lem:1}.
\qed

\begin{theorem}
\label{cor:1} Let $V$ denote a finite-dimensional irreducible
$\boxtimes$-module.  Let $r,s,t$ denote mutually distinct elements
in $\I$. Then for all $\lambda \in \Bbb K$ we have
\begin{eqnarray}
V_{rs}(\lambda) + V_{rs}(\lambda-2) + \cdots = V_{rt}(\lambda) +
V_{rt}(\lambda-2) + \cdots . \label{eq:eq}
\end{eqnarray}
We are using Notation \ref{note:u}.
\end{theorem}

{\it Proof.} Let $Y$ denote the left-hand side of (\ref{eq:eq})
and let $Z$ denote the right-hand side of (\ref{eq:eq}).  We show
$Y = Z$.  We first show $Y \subseteq Z$. Since the dimension of
$V$ is finite, there exists a nonnegative integer $j$ such that
$V_{rs}(\lambda-2i) = 0$ and $V_{rt}(\lambda-2i) = 0$ for all
integers $i
> j$. Observe that $Y= \sum_{i=0}^j V_{rs}(\lambda -2i)$ and $Z=
\sum_{i=0}^{j} V_{rt}(\lambda -2i)$.  By (\ref{eq:es2}) we find
$Z$ is the set of vectors in $V$ on which
\begin{eqnarray}
\prod_{i=0}^{j} (X_{rt} - (\lambda-2i)I) \label{eq:in}
\end{eqnarray}
vanishes.  Using the table in Theorem \ref{lem:1} (row $t=r,u \not=
s$) we find that (\ref{eq:in}) vanishes on
$V_{rs}(\lambda-2i)$ for $0 \leq i \leq j$.  Therefore
$V_{rs}(\lambda-2i) \subseteq Z$ for $0 \leq i \leq j$ so $Y
\subseteq Z$.

To get $Z \subseteq Y$ interchange the roles of $s,t$ in the
argument so far.  We conclude that $Y=Z$ and the result follows.
\qed

\begin{corollary}
\label{cor:2} Let $V$ denote a finite-dimensional irreducible
$\boxtimes$-module and fix $\lambda \in \Bbb K$.  Then for
distinct $ r,s \in  \I$ the dimension of $V_{rs}(\lambda)$ is
independent of $r,s$. We are using Notation \ref{note:u}.
\end{corollary}

{\it Proof.}  We define a binary relation on the generators of
$\boxtimes$ called {\it adjacent}.  By definition two distinct
generators $X_{rs},X_{tu}$ are adjacent whenever $r=t$ or $s=u$.
We observe the adjacency relation is symmetric.  We further
observe that the generators of $\boxtimes$ are connected with
respect to adjacency. By these comments it suffices to show
$V_{rs}(\lambda), V_{tu}(\lambda)$ have the same dimension for
adjacent generators $X_{rs},X_{tu}$. First assume $r=t$.  By
Theorem \ref{cor:1} we find
\begin{eqnarray}
V_{rs}(\lambda) + V_{rs}(\lambda-2) + \cdots = V_{tu}(\lambda) +
V_{tu}(\lambda-2) + \cdots . \label{eq:star2}
\end{eqnarray}
Applying Theorem \ref{cor:1} (with
$\lambda$ replaced by $\lambda -2$), we find
\begin{eqnarray}
V_{rs}(\lambda-2) + V_{rs}(\lambda-4) + \cdots =  V_{tu}(\lambda-2)
+V_{tu}(\lambda-4) + \cdots .  \label{eq:star}
\end{eqnarray}
Let $H$ denote the sum on either side of (\ref{eq:star}).  Comparing
(\ref{eq:star2}) and (\ref{eq:star}) we find
\begin{eqnarray}
V_{rs}(\lambda) + H = V_{tu}(\lambda) + H. \label{eq:star1}
\end{eqnarray}
The sum on either side of (\ref{eq:star1}) is direct, so
$V_{rs}(\lambda)$ and $V_{tu}(\lambda)$ have the same dimension.

Next we assume $s=u$.  Applying what we have done so
far with $(r,s,t,u,\lambda)$ replaced by $(s,r,u,t, -\lambda)$ we find
$V_{sr}(-\lambda)$
and $V_{ut}(-\lambda)$ have the same dimension.  Recall that $X_{sr} =
-X_{rs}$ and $X_{ut} = -X_{tu}$ so $V_{sr}(-\lambda) =
V_{rs}(\lambda)$ and $V_{ut}(-\lambda) =
V_{tu}(\lambda)$.  By these comments $V_{rs}(\lambda)$ and $V_{tu}(\lambda)$
have the same dimension.  We have now shown $V_{rs}(\lambda),
V_{tu}(\lambda)$ have the
same dimension for adjacent generators $X_{rs},X_{tu}$ and the result
follows.
\qed

\begin{corollary}
\label{cor:samedim} Let $V$ denote a finite-dimensional
irreducible $\boxtimes$-module and fix $\lambda \in \Bbb K$.  Then
for distinct $r,s \in \I$ the spaces $V_{rs}(\lambda),
V_{rs}(-\lambda)$ have the same dimension.  We are using Notation
\ref{note:u}.
\end{corollary}

{\it Proof.}  By Corollary \ref{cor:2} we find
$V_{rs}(\lambda),V_{sr}(\lambda)$ have the same dimension.  But
$X_{rs} = -X_{sr}$ so $V_{sr}(\lambda) = V_{rs}(-\lambda)$. \qed

\begin{theorem}
\label{lem:2facts}
Let $V$
denote a finite-dimensional irreducible $\boxtimes$-module.  Then
the following $(i),(ii)$ hold.
\begin{enumerate}
\item
For distinct $r,s \in \I$ the generator $X_{rs}$ is diagonalizable
on $V$.
\item
There exists an integer $d \geq 0$ such that for all distinct $r,s
\in \I$ the eigenvalues for $X_{rs}$ on $V$ are $d, d-2, \ldots,
-d$.
\end{enumerate}
\end{theorem}

{\it Proof.} Throughout this proof we will use Notation
\ref{note:u}. Let $r,s$ denote distinct elements of $\I$.  Since
$\Bbb K$ is algebraically closed, all the eigenvalues for $X_{rs}$
on $V$ are contained in $\Bbb K$.  Since $V$ has positive finite
dimension, the action of $X_{rs}$ on $V$ has at least one
eigenvalue $\lambda$.  Since ${\rm Char} (\Bbb K) = 0$, the
scalars $\lambda, \lambda +2, \ldots$ are mutually
distinct, so they cannot all be eigenvalues for $X_{rs}$ on $V$;
consequently there exists $\theta \in \Bbb K$ such that
$V_{rs}(\theta) \not= 0$ but $V_{rs}(\theta+2) = 0$.  Similarly
the scalars  $\theta, \theta -2, \ldots$ are mutually
distinct, so they cannot all be eigenvalues for $X_{rs}$ on $V$;
consequently there exists an integer $d \geq 0$ such that
$V_{rs}(\theta - 2i)$ is nonzero for $0 \leq i \leq d$ and zero
for $i=d+1$.  We will now show that
\begin{eqnarray}
V_{rs}(\theta) + \cdots + V_{rs}(\theta -
2d) \label{eq:vis}
\end{eqnarray}
is equal to $V$.

By Corollary \ref{cor:act} and since $V_{rs}(\theta+2) = 0$,
$V_{rs}(\theta - 2d - 2) = 0$ we find that (\ref{eq:vis}) is
$X_{tu}$-invariant for all distinct $t,u \in \I$.  Therefore
(\ref{eq:vis}) is $\boxtimes$-invariant.  Recall the
$\boxtimes$-module $V$ is irreducible so (\ref{eq:vis}) is equal
to either $0$ or $V$. By construction each term in (\ref{eq:vis})
is nonzero and there is at least one term, so (\ref{eq:vis}) is
nonzero.  Therefore (\ref{eq:vis}) is equal to $V$.  This shows
that the action of $X_{rs}$ on $V$ is diagonalizable with
eigenvalues $\Delta = \{ \theta -2i | 0 \leq i \leq d \} $.

It remains to show $\theta = d$.  By Corollary \ref{cor:samedim} we
find $\Delta = - \Delta$.  It follows that $\theta = -\theta + 2d$ so
$\theta = d$.
\qed

\begin{definition}
\rm Let $V$ denote a finite-dimensional irreducible
$\boxtimes$-module. By the {\it diameter} of $V$ we mean the
nonnegative integer $d$ from Theorem \ref{lem:2facts}$(ii)$.
\end{definition}

Our discussion of $\boxtimes$-modules will continue after we recall
the notion of a flag.

\section{Flags}
\label{sect:flag}

\begin{definition}
\rm Let $V$ denote a vector space over $\Bbb K$ with positive
finite dimension.  Let $d$ denote a nonnegative integer.  By a
{\it flag on $V$ of diameter $d$}, we mean a sequence $F_0,F_1,
\ldots, F_d$ consisting of mutually distinct subspaces of $V$ such
that $F_0 \not= 0$, $F_{i-1} \subseteq F_i$ for $1 \leq i \leq d$,
and $F_d = V$.  We call $F_i$ the $i^{th}$ {\it component} of the
flag.
\end{definition}

Let $V$ denote a vector space over $\Bbb K$
with positive finite dimension.  Let $d$ denote a
nonnegative integer.  By a {\it decomposition of $V$ of diameter $d$}, we
mean a sequence $V_0,V_1, \ldots, V_d$ consisting of nonzero subspaces
of $V$ such that
$$V = V_0 + \cdots + V_d \qquad \mbox{(direct sum).}$$
We do not assume each of $V_0, \ldots, V_d$ has dimension $1$.
For $0 \leq i \leq d$ we call $V_i$ the {\it $i^{th}$ subspace} of the
decomposition.

The following construction yields a flag on $V$.  Let $V_0,
\ldots, V_d$ denote a decomposition of $V$.  Set
$$F_i = V_0 + \cdots + V_i$$
for $0 \leq i \leq d$.  Then the sequence $F_0, \ldots, F_d$ is a
flag on $V$.  We say this flag is {\it induced} by $V_0, \ldots, V_d$.

Let $V_0, \ldots, V_d$ denote a decomposition of $V$.  By the {\it
inversion} of this decomposition, we mean the decomposition $V_d,
\ldots, V_0$ of $V$.

We now discuss the notion of {\it opposite} flags.  Let $F$ and $G$
denote flags on $V$.  These flags are said to be {\it opposite} whenever
there exists a decomposition $V_0, \ldots, V_d$ of $V$ such that $F$ is
induced by $V_0, \ldots, V_d$ and $G$ is induced by $V_d,
\ldots, V_0$. In this case
\begin{eqnarray}
\label{eq:note1}
F_i \cap G_j = 0 \qquad \mbox{if } i+j < d
\end{eqnarray}
and
\begin{eqnarray}
\label{eq:note2}
V_i = F_i \cap G_{d-i} \qquad \mbox{for } 0 \leq i \leq d.
\end{eqnarray}
In particular $V_0, \ldots, V_d$ is uniquely determined by the
ordered pair $F,G$.  We say $V_0, \ldots, V_d$ is {\it
induced} by $F,G$.

\section{Finite-Dimensional Irreducible $\boxtimes$-Modules,
  Revisited}

We return our attention to a finite-dimensional irreducible
$\boxtimes$-module $V$.  Recall the set $\I$ from Definition
\ref{def:main} has four elements.  With each element of $\I$ we
associate a flag on $V$.  We start with a definition.

\begin{definition}
\rm Let $V$ denote a finite-dimensional irreducible
$\boxtimes$-module.  Let $(r,s)$ denote an ordered pair of
distinct elements of $\I$.  By Theorem \ref{lem:2facts} the
sequence $V_{rs}(-d),V_{rs}(2-d), \ldots, V_{rs}(d)$ is a
decomposition of $V$, where $d$ denotes the diameter.  We will
call this the {\it decomposition of $V$ associated with $(r,s)$}.
\end{definition}

\begin{note}
\label{note:29}
\rm For distinct $r,s \in \I$ the decomposition of $V$ associated
with $(s,r)$ is the inversion of the decomposition of $V$
associated with $(r,s)$.
\end{note}

\begin{lemma}
\label{lem:nwaa} Let $V$ denote a finite-dimensional irreducible
$\boxtimes$-module. Let $r,s$ denote distinct elements of $\I$ and
consider the decomposition of $V$ associated with $(r,s)$.  Then
the flag induced by this decomposition is independent of $s$.
\end{lemma}

{\it Proof.}  Let $t$ denote an element of $\I$ such that $r,s,t$
are mutually distinct.  It suffices to show that the flag induced
by the decomposition of $V$ associated with $(r,s)$ equals the
flag induced by the decomposition of $V$ associated with $(r,t)$.
Let $d$ denote the diameter of $V$.  For $0 \leq i \leq d$ the
$i^{th}$ component of the first flag is $V_{rs}(-d)
+ \cdots + V_{rs}(2i-d)$ and the $i^{th}$ component of the second
flag is $V_{rt}(-d) + \cdots + V_{rt}(2i-d)$.  These
components are equal by Theorem \ref{cor:1} (with $\lambda =
2i-d$). Therefore the flags are equal. \qed

\medskip

\begin{definition}
\label{def:nwaa} \rm Let $V$ denote a finite-dimensional
irreducible $\boxtimes$-module. For $r \in \I$, by the {\it flag on
$V$ associated with $r$} we mean the flag discussed in Lemma
\ref{lem:nwaa}.
\end{definition}

\medskip

The next corollary restates Lemma \ref{lem:nwaa} in light of
Definition \ref{def:nwaa}.

\begin{corollary}
\label{cor:nwbb} Let $V$ denote a finite-dimensional irreducible
$\boxtimes$-module.  For $r \in \I$ consider the flag on $V$ associated
with $r$.  The components of this flag are described as follows.
Let $d$ denote the diameter of $V$ and pick $s \in \I$ such that
$r \not= s$. Then for $0 \leq i \leq d$ the $i^{th}$ component of
the flag is
\begin{eqnarray}
\label{eq:nwq} V_{rs}(-d)+V_{rs}(2-d) + \cdots + V_{rs}(2i-d).
\end{eqnarray}
\end{corollary}

\begin{lemma}
\label{lem:nwcc} Let $V$ denote a finite-dimensional irreducible
$\boxtimes$-module.  Recall from Definition \ref{def:main} that
$\I$ has four elements and consider the corresponding flags on $V$
from Definition \ref{def:nwaa}.  These four flags are mutually
opposite.
\end{lemma}

{\it Proof.}  For distinct $r,s \in \I$ let
$F$ denote the flag on $V$ associated with $r$ and let $G$ denote
the flag on $V$ associated with $s$.  We show $F$ and $G$ are
opposite.  By Note \ref{note:29} the decomposition of $V$ associated
with $(r,s)$ is the inversion of the decomposition of $V$ associated
with $(s,r)$.  The first decomposition induces $F$ and the second
decomposition induces $G$.  It follows that $F$ and $G$ are
opposite. \qed

\begin{lemma}
\label{lem:nwdd} Let $V$ denote a finite-dimensional irreducible
$\boxtimes$-module of diameter $d$.  Let $r,s$ denote distinct elements
of $\I$.  For $0 \leq i \leq d$ the subspace $V_{rs}(2i-d)$
is the intersection of the following two sets:
\begin{enumerate}
\item The $i^{th}$ component of the flag on $V$ associated with $r$.
\item The $(d-i)^{th}$ component of the flag on $V$ associated with $s$.
\end{enumerate}
\end{lemma}

{\it Proof.}  By Definition \ref{def:nwaa} the flag on $V$
associated with $r$ is induced by $V_{rs}(-d),
\ldots, V_{rs}(d)$ and the flag on $V$ associated with $s$ is
induced by $V_{rs}(d), \ldots, V_{rs}(-d)$.  The
result follows. \qed

\medskip

\section{From $\boxtimes$-Modules to $\O$-Modules}

In this section we prove Theorem \ref{th:H1}.  We begin with a
lemma.

\begin{lemma}
\label{lem:ne} Let $V$ denote a finite-dimensional irreducible
$\boxtimes$-module. Let $r,s,t,u$ denote mutually distinct
elements of $\I$. Let $W$ denote a nonzero subspace of $V$ such
that $X_{rs}W \subseteq W$, $X_{tu}W \subseteq W$.  Then $W = V$.
\end{lemma}

{\it Proof.}  Without loss of generality, we assume $W$ is irreducible
as a module for $X_{rs},X_{tu}$.  To show $W = V$, it suffices to show
that $W$ is $\boxtimes$-invariant.  To this end, we show
\begin{eqnarray}
\label{eq:cl} X_{ab} W \subseteq W
\end{eqnarray}
for all distinct $a,b \in \I$.  Assume $ab$ is not one of $rs, tu,
sr, ut$; otherwise (\ref{eq:cl}) holds by the construction and
(\ref{eq:1}).  Replacing an appropriate subset of $\{ rs,tu,ab \}$
by the corresponding subset of $\{ sr,ut,ba \}$ we may assume,
without loss of generality, that $a=s$ and $b=t$. Therefore it
suffices to show $X_{st} W \subseteq W$. We define $W' :=\{ w \in
W | X_{st} w \in W \}$ and show $W' = W$. By (\ref{eq:2}) we
routinely find that $X_{rs} W' \subseteq W'$ and $X_{tu} W'
\subseteq W'$. By these comments and the irreducibility of $W$ we
find either $W'=0$ or $W'=W$.

\noindent
{\it Claim.} $W' \not= 0$.

\noindent {\it Proof of claim.}  Let $d$ denote the diameter of $V$.
Define $W_j = W \cap
V_{rs}(2j-d)$ for $0 \leq j \leq d$. The nonzero spaces among $W_0
, \ldots , W_d$ are eigenspaces of $X_{rs}$ on $W$. By Theorem
\ref{lem:2facts} we find $W = \sum_{j=0}^d W_j$.   By this and
since $W \not= 0$ we find $W_0 , \ldots , W_d$ are not all $0$.
Define $ m = \mbox{min} \{ j | 0 \leq j \leq d, W_j \not= 0 \}.$

Define $W_j^* = W \cap V_{tu}(2j-d)$ for $0 \leq j \leq d$.  The
nonzero spaces among $W_0^* , \ldots , W_d^*$ are eigenspaces of
$X_{tu}$ on $W$. By Theorem \ref{lem:2facts} we find $W =
\sum_{j=0}^d W_j^*$.  By this and since $W \not= 0$ we find $W_0^*
, \ldots , W_d^*$ are not all $0$. Define $ n = \mbox{min} \{ j | 0 \leq
j \leq d, W_j^* \not= 0 \}. $

We show $m=n$.  To do this we first assume $m < n$ and get a
contradiction. By construction $W_j^* \subseteq V_{tu}(2j-d)$ for
$0 \leq j \leq d$ and $W = \sum_{j=n}^d W_j^*$ so $W \subseteq
\sum_{j=n}^d V_{tu}(2j-d).$ Therefore $W$ is contained in the
$(d-n)^{th}$ component of the flag on $V$ associated with $u$. By
construction $W_m \subseteq V_{rs}(2m-d)$ so $W_m \subseteq
\sum_{j=0}^m V_{rs}(2j-d)$.  Therefore $W_m$ is contained in the
$m^{th}$ component of the flag on $V$ associated with $r$. Recall
by Lemma \ref{lem:nwcc} that the flag on $V$ associated with $u$
and the flag on $V$ associated with $r$ are opposite. By
(\ref{eq:note1}) the $(d-n)^{th}$ component of the flag on $V$
associated with $u$ and the $m^{th}$ component of the flag on $V$
associated with $r$ have $0$ intersection. Therefore $W \cap W_m =
0$, contradicting the fact that $W_m$ is nonzero in $W$. Therefore
$m \geq n$.

Next we assume $m > n$ and get a contradiction.  By construction
$W_j \subseteq V_{rs}(2j-d)$ for $0 \leq j \leq d$ and $W =
\sum_{j=m}^d W_j$ so $W \subseteq \sum_{j=m}^d V_{rs} (2j-d)$.
Therefore $W$ is contained in the $(d-m)^{th}$ component of the
flag on $V$ associated with $s$. By construction $W_n^* \subseteq
V_{tu}(2n-d)$ so $W_n^* \subseteq \sum_{j=0}^n V_{tu}(2j-d)$.
Therefore $W_n^*$ is contained in the $n^{th}$ component of the
flag on $V$ associated with $t$. Recall by Lemma \ref{lem:nwcc}
that the flag on $V$ associated with $s$ and the flag on $V$
associated with $t$ are opposite.  By (\ref{eq:note1}) the
$(d-m)^{th}$ component of the flag on $V$ associated with $s$ and
the $n^{th}$ component of the flag on $V$ associated with $t$ have
$0$ intersection. Therefore $W \cap W_n^* = 0$, contradicting the
fact that $W_n^*$ is nonzero in $W$.  Therefore $m \leq n$.

Combining the previous two paragraphs we find $m = n$.  By the
previous paragraph we find $W_n^*$ is contained in the
intersection of the $(d-m)^{th}$ component of the flag on $V$
associated with $s$ and the $m^{th}$ component of the flag on $V$
associated with $t$.  By Lemma \ref{lem:nwdd} we find $W_n^*
\subseteq V_{st}(d-2m)$.  Since $V_{st}(d-2m)$ is an eigenspace
for $X_{st}$ we find $X_{st} W_n^* \subseteq W_n^*$. Therefore
$W_n^* \subseteq W'$, so $W' \not= 0$ and the claim is proved.

This shows that $W' = W$ and therefore $X_{st}W \subseteq W$.  By this
and the comments prior to the claim we have now shown that
$X_{ab} W \subseteq W$ for all distinct $a,b \in \I$.
Therefore $W$ is $\boxtimes$-invariant.  Since the $\boxtimes$-module
$V$ is irreducible we have $W=V$. \qed

\medskip

We are now ready to prove Theorem \ref{th:H1}.

\medskip

{\it Proof of Theorem \ref{th:H1}.} Since $X_{01},X_{23}$ satisfy
(\ref{eq:fo1}) and (\ref{eq:fo2}) there exists an $\O$-module
structure on $V$ such that $X,Y$ act on $V$ as
$X_{01},X_{23}$ respectively.  This structure is unique since $X,Y$
generate $\O$.  This structure is irreducible by
Lemma \ref{lem:ne}.  This structure has type $(0,0)$ by Theorem
\ref{lem:2facts}$(ii)$. \qed

\section{From $\O$-Modules to $\boxtimes$-Modules}

In this section we prove Theorem \ref{th:H2}.  We refer to the
following setup.

\begin{definition}
\rm \label{def:asm} Let $V$ denote a finite-dimensional
irreducible $\O$-module of type $(0,0)$.  Let $A,A^*$ denote the
actions on $V$ of the generators $X,Y$ respectively and recall by
Corollary \ref{cor:t0o} that $A,A^*$ is a tridiagonal pair on $V$.
Let $\I$ denote the set consisting of the four symbols $\{0,1,2,3
\}$. With each element of $\I$ we will associate a flag on $V$ of
diameter $d$ where $d$ denotes the diameter of $V$.  For $0 \leq i
\leq d$ the $i^{th}$ component of this flag is given as follows.

\medskip

\begin{tabular}{c|c}
Element of $\I$ & $i^{th}$ component of associated flag \\
\hline
$0$ & $V_A(-d)+ \cdots + V_A(2i-d)$ \\
$1$ & $V_A(d)+ \cdots + V_A(d-2i)$ \\
$2$ &  $V_{A^*}(-d)+  \cdots + V_{A^*}(2i-d)$ \\
$3$ & $V_{A^*}(d)+  \cdots + V_{A^*}(d-2i)$
\end{tabular}
\end{definition}

\begin{lemma}
\label{lem:asm}
\cite[Theorem 7.3]{ter3} The four flags from Definition
\ref{def:asm} are mutually opposite.
\end{lemma}

\begin{definition}
\rm \label{def:12} Adopt the assumptions of Definition
\ref{def:asm}. For distinct $r,s \in \I$ let $x_{rs}:V \to V$
denote the unique linear transformation satisfying the following
conditions.  By Definition \ref{def:asm} each of $r,s$ is
associated with a flag on $V$ of diameter $d$.  By Lemma \ref{lem:asm}
the flags $r,s$ are opposite; let $ U_0, \ldots, U_d$ denote the
induced decomposition of $V$.  For $0 \leq i \leq d$,
$U_i$ is an eigenspace for $x_{rs}$ with eigenvalue $2i-d$.
\end{definition}

We will be using the following notation.

\begin{notation}
\rm \label{not:11} With reference to Definition \ref{def:12}, for
distinct $r,s \in \I$ and for all $\lambda \in \Bbb K$ we define
$$V_{rs}(\lambda) = \{ v \in V | x_{rs} v = \lambda v \}.$$
\end{notation}

The next two lemmas illustrate the notation in Definition
\ref{def:12} and will be useful later in the paper.

\begin{lemma}
\label{lem:note} With reference to Definition \ref{def:12}, we
have $x_{01} = A$ and $x_{23} = A^*$.
\end{lemma}

{\it Proof.}  Immediate from Definition \ref{def:12}.
\qed

\begin{lemma}
\label{lem:go} With reference to Definition \ref{def:12}, for $r
\in \I$ consider the flag on $V$ associated with $r$. The
components of this flag are described as follows. Pick $s \in \I$
such that $r \not= s$. Then for $0 \leq i \leq d$ the $i^{th}$
component of the flag is
$$V_{rs}(-d)+ \cdots + V_{rs}(2i-d).$$
We are using Notation \ref{not:11}.
\end{lemma}

{\it Proof.}  By Definition \ref{def:12} the decomposition
$V_{rs}(-d), \ldots, V_{rs}(d)$ of $V$ is induced by
the flag on $V$ associated with $r$ and the flag on $V$ associated
with $s$. By the discussion at the end of Section \ref{sect:flag}
the flag $r$ is induced by $V_{rs}(-d), \ldots,
V_{rs}(d)$. The result follows. \qed

\medskip

We will need the following three lemmas.

\begin{lemma}
\label{note:minus} With reference to Definition \ref{def:12}, for
distinct $r,s \in \I$ we have 
$$ x_{rs} + x_{sr} =0. $$
\end{lemma}

{\it Proof.}  Let $U_0, \ldots , U_d$ denote the decomposition
induced by the flag on $V$ associated with $r$ and the flag on $V$
associated with $s$.  By Definition \ref{def:12} and with
reference to Notation \ref{not:11}, we find $V_{rs}(2i-d) = U_i =
V_{sr}(d-2i)$ for $0 \leq i \leq d$. The result follows. \qed

\begin{lemma}
\label{lem:ii} With reference to Definition \ref{def:12}, for mutually
distinct $r,s,t \in \I$ we have 
\begin{eqnarray}
[x_{rs},x_{st}] = 2x_{rs} + 2x_{st}. \label{eq:r}
\end{eqnarray}
\end{lemma}

{\it Proof.}  Throughout this proof we will use Notation \ref{not:11}.

We break our proof into the following three cases:

(i) $st$ is one of $01, 10, 23, 32$.

(ii) $rs$ is one of $01, 10, 23, 32$.

(iii) Neither $rs$ nor $st$ is one of $01,10,23,32$.

Case (i):  We will invoke Lemma \ref{lem:conf2}.  We begin with some
comments.  Recall that
$x_{rs}$ is diagonalizable on $V$ with eigenvalues $d, d-2,
\ldots, -d$. We now show that
\begin{eqnarray}
\label{eq:newe1} (x_{st} + (2i-d)I)V_{rs}(2i-d) \subseteq
V_{rs}(2i-d+2)
\end{eqnarray}
for $0 \leq i \leq d$.  Let $u$ denote the unique element of $\I$
such that $r,s,t,u$ are mutually distinct. By Lemma \ref{lem:note}
and Lemma \ref{note:minus} we find that $x_{st},x_{ur}$ is either
$\epsilon A, \epsilon^* A^*$ or $\epsilon^* A^*,\epsilon A$ for
some $\epsilon, \epsilon^* \in \{ 1,-1 \}$.  By this and Definition
\ref{def:asm} we find $x_{st},x_{ur}$
is a tridiagonal pair on $V$ of diameter $d$ and $d, d-2, \ldots ,
-d$ is both an eigenvalue sequence and dual eigenvalue sequence of
$x_{st},x_{ur}$. By Lemma \ref{lem:go} and using \cite[Theorem
4.6]{ITT} we find that (\ref{eq:newe1}) holds for $0 \leq i \leq
d$.  By this and Lemma \ref{lem:conf2} we find (\ref{eq:r}) holds.

Case (ii):  Observe that $sr$ is one of $01, 10, 23, 32$.
By case (i) of the current proof we find that
\begin{eqnarray}
\label{eq:fy} [x_{ts},x_{sr}] = 2x_{ts}+2x_{sr}.
\end{eqnarray}
By Lemma \ref{note:minus} we find $x_{ts} = -x_{st}$ and $x_{sr}
= -x_{rs}$.  Substituting this into (\ref{eq:fy}) and simplifying
the result we find that (\ref{eq:r}) holds.

Case (iii):  Similarly to case (i) we will invoke Lemma
\ref{lem:conf2}.  Recall that
$x_{rs}$ is diagonalizable on $V$ with eigenvalues $d, d-2,
\ldots, -d$. We now show that
\begin{eqnarray}
\label{eq:newe2} (x_{st} + (2i-d)I) V_{rs}(2i-d) \subseteq
V_{rs}(2i-d+2)
\end{eqnarray}
for $0 \leq i \leq d$.  Observe $rt$ is one of $01, 10, 23, 32$. By
case (ii) of the current proof we find that $[x_{rt},x_{ts}] = 2x_{rt} +
2x_{ts}$.  By Lemma \ref{lem:conf2} (with $\lambda = d-2i$) and since
$x_{rt} = -x_{tr}$, $x_{ts}= -x_{st}$ we find 
\begin{eqnarray}
(x_{st} + (2i-d) I)V_{tr}(2i-d) \subseteq V_{tr}(2i-d -2) \qquad 0
\leq i \leq d. \label{eq:d3}
\end{eqnarray}
For $0 \leq i \leq d$ we have
\beast (x_{st} + (2i-d)
I)V_{rs}(2i-d) & \subseteq & (x_{st} + (2i-d) I)
\sum_{n=0}^i V_{rs}(2n-d) \\
& = & (x_{st} + (2i-d) I) \sum_{n=0}^i V_{tr}(d-2n) \qquad \qquad
\mbox{by Lemma \ref{lem:go}} \\
& \subseteq & \sum_{n=0}^{i+1} V_{tr}(d-2n)
\qquad \qquad \qquad \qquad \qquad \quad \mbox{by (\ref{eq:d3})} \\
& = & \sum_{n=0}^{i+1}V_{rs}(2n-d) \qquad \qquad \qquad \qquad \qquad
\quad \mbox{by Lemma
\ref{lem:go}} \eeast
and
\beast (x_{st} + (2i-d) I)V_{rs}(2i-d) &
\subseteq & (x_{st} + (2i-d) I)
\sum_{n=i}^dV_{rs}(2n-d) \\
& = & (x_{st} + (2i-d) I) \sum_{n=i}^d  V_{st}(d-2n) \qquad \qquad
\mbox{by Lemma \ref{lem:go}}\\
& = & \sum_{n=i+1}^d  V_{st}(d-2n) \\
& = & \sum_{n=i+1}^d V_{rs}(2n-d) \qquad \qquad \qquad \qquad \qquad
\mbox{by Lemma \ref{lem:go}}. \eeast

Combining these observations we find that (\ref{eq:newe2}) holds
for $0 \leq i \leq d$.  By this and Lemma \ref{lem:conf2} we find
that (\ref{eq:r}) holds. \qed

\begin{lemma}
\label{lem:iii} With reference to Definition \ref{def:12}, for
mutually distinct $r,s,t,u \in \I$ we have
\begin{eqnarray}
[x_{rs},[x_{rs},[x_{rs},x_{tu}]]] = 4[x_{rs},x_{tu}]. \label{eq:s}
\end{eqnarray}
\end{lemma}

{\it Proof.}   Throughout this proof we will use Notation \ref{not:11}.

We will invoke Lemma \ref{lem:conf}.  We begin with some comments.
Recall that $x_{rs}$ is diagonalizable on $V$ with eigenvalues $d,
d-2, \ldots, -d$.  We now show that
\begin{eqnarray}
\label{eq:newe3} x_{tu}V_{rs}(2i-d) \subseteq V_{rs}(2i-d+2)+
V_{rs}(2i-d)+ V_{rs}(2i-d-2)
\end{eqnarray}
for $0 \leq i \leq d$.  By Lemma \ref{lem:ii} we find $[x_{st},x_{tu}]
= 2x_{st} + 2x_{tu}$. By this and Lemma \ref{lem:conf2} we find
\begin{eqnarray}
(x_{tu} + (2i-d) I)V_{st}(2i-d) \subseteq  V_{st}(2i-d +2)
\label{eq:e4} \qquad 0 \leq i \leq d.
\end{eqnarray}
By Lemma \ref{lem:ii} we find $[x_{ru},x_{ut}] =
2x_{ru}+2x_{ut}$.  By Lemma \ref{lem:conf2} (with $\lambda = d-2i$)
and since $x_{ru} = -x_{ur}$, $x_{ut} = -x_{tu}$ we find 
\begin{eqnarray}
(x_{tu} + (2i-d) I)V_{ur}(2i-d) \subseteq V_{ur}(2i-d -2) \qquad 0
\leq i \leq d. \label{eq:e3}
\end{eqnarray}
For $0 \leq i \leq d$ we have
\beast
x_{tu}V_{rs}(2i-d) & \subseteq & x_{tu} \sum_{n=i}^d V_{rs}(2n-d) \\
& = & x_{tu} \sum_{n=i}^d V_{st}(d-2n)
\qquad \qquad \mbox{by Lemma \ref{lem:go}}\\
& \subseteq &  \sum_{n=i-1}^dV_{st}(d-2n)
\qquad \qquad \quad \mbox{by (\ref{eq:e4})} \\
 & = & \sum_{n=i-1}^dV_{rs}(2n-d) \qquad \qquad \quad
 \mbox{by Lemma \ref{lem:go}}
\eeast
and
\beast
x_{tu} V_{rs}(2i-d) & \subseteq & x_{tu} \sum_{n=0}^i V_{rs}(2n-d) \\
& = & x_{tu} \sum_{n=0}^i V_{ur}(d-2n)
\qquad \qquad \mbox{by Lemma \ref{lem:go}}\\
& \subseteq & \sum_{n=0}^{i+1}V_{ur}(d-2n)
\qquad \qquad \quad  \mbox{by (\ref{eq:e3})} \\
 & = & \sum_{n=0}^{i+1}  V_{rs}(2n-d) \qquad
 \qquad \quad  \mbox{by Lemma \ref{lem:go}}.
\eeast Combining these observations we find that (\ref{eq:newe3})
holds for all $0 \leq i \leq d$.   By this and Lemma
\ref{lem:conf} we find that (\ref{eq:s}) holds. \qed

\medskip

\begin{theorem}
\label{th:pth} With reference to Definition \ref{def:12}, there
exists a unique $\boxtimes$-module structure on $V$ such that
the generator $X_{rs}$ acts on $V$ as $x_{rs}$ for all distinct $r,s
\in \I$. This
$\boxtimes$-module structure is irreducible.
\end{theorem}

{\it Proof.}  Using Lemma \ref{note:minus}, Lemma
\ref{lem:ii}, and Lemma \ref{lem:iii} we find that there exists a
$\boxtimes$-module structure on $V$ such that $X_{rs}$ acts on $V$
as $x_{rs}$ for distinct $r,s \in \I$.  This module structure is
unique since the set $\{ X_{rs} | r,s \in \I, r \not= s \}$ is a
generating set for $\boxtimes$.  By Lemma \ref{lem:note} we find
that $A,A^*$ are among the linear transformations $x_{rs}$, $r,s
\in \I$, $r \not= s$. Recall that $A,A^*$ is a tridiagonal pair on
$V$. Therefore the $\boxtimes$-module $V$ is irreducible by Definition
\ref{def:tp}(iv).
\qed

\medskip

We are now ready to prove Theorem \ref{th:H2}.

\medskip

{\it Proof of Theorem \ref{th:H2}.}  Consider the $\boxtimes$-module
structure on $V$ from Theorem \ref{th:pth}.  In this module structure
$X_{01},X_{23}$ act on $V$ as $X,Y$ respectively by Definition
\ref{def:asm} and Lemma \ref{lem:note}.  Next we show that this
$\boxtimes$-module structure is unique. Suppose we are given any
$\boxtimes$-module structure on $V$ where $X_{01},X_{23}$ act on $V$
as $X,Y$ respectively.  This module structure is irreducible since the
$\O$-module $V$ is irreducible.  For each generator
$X_{rs}$ of $\boxtimes$ the action on $V$ is determined by the
decomposition of $V$ associated with $(r,s)$.  By Lemma
\ref{lem:nwdd} the decomposition of $V$ associated with $(r,s)$ is
determined by the flag on $V$ associated with $r$ and the flag on
$V$ associated with $s$. Therefore our $\boxtimes$-module
structure on $V$ is determined by the four flags associated with
the four elements of $\I$.  By Corollary \ref{cor:nwbb} the four
flags associated with the four elements of $\I$ are determined by
the actions of $X_{01}$ and $X_{23}$ on $V$.  Therefore the given
$\boxtimes$-module structure on $V$ is determined by the actions
of $X_{01}$ and $X_{23}$ on $V$, so the $\boxtimes$-module
structure is unique.  We have now shown there exists a unique
$\boxtimes$-module structure on $V$ where $X_{01},X_{23}$ act on
$V$ as $X,Y$ respectively.  We mentioned earlier that this module
structure is irreducible. \qed

\bigskip

\noindent
{\bf Acknowledgements:}  This paper was written while the author was
a graduate student at the University of Wisconsin-Madison.  The author
would like to thank his advisor Paul Terwilliger for his many valuable
ideas and suggestions.

\noindent Brian Hartwig  \hfil \break Department of Mathematics
\hfil \break University of Wisconsin  \hfil \break 480 Lincoln
Drive \hfil \break Madison, Wisconsin, 53706, USA \hfil \break
email: hartwig@math.wisc.edu \hfil \break


\begin{thebibliography}{9}

\bibitem{A}
H. ~Au-Yang, B.M. ~McCoy, J.H.H. ~Perk, S. ~Tang, and M.L. Yan.
\newblock Commuting transfer matrices in chiral Potts modules:
Solutions of the star-triangle equations with genus greater than $1$.
\newblock {\it Phys. Lett. A} {\bf 123} (1987) 219-223.

\bibitem{bandt}
G. ~Benkart and P. ~Terwilliger.
\newblock Irreducible modules for the quantum affine algebra $U_q
(\widehat{sl}_2)$ and its Borel subalgebra.
\newblock {\it J. Algebra} {\bf 282} (2004), no. 1, 172--194.

\bibitem{DR}
E. ~Date and S.S ~Roan.
\newblock The structure of quotients of the Onsager algebra by
closed ideals.
\newblock {\it J. Phys. A: Math. Gen.} {\bf 33} (2000) 3275--3296.

\bibitem{Dav1}
B. ~Davies.
\newblock Onsager's algebra and superintegrability.
\newblock {\it J. Phys. A:Math. Gen.} {\bf 23} (1990), 2245--2261.

\bibitem{Dav2}
B. ~Davies.
\newblock Onsager's algebra and the Dolan-Grady condition in the
non-self-dual case.
\newblock {\it J. Math. Phys.} {\bf 32} (1991), 2945--2950.

\bibitem{DG}
L. ~Dolan and M. ~Grady.
\newblock Conserved charges from self-duality.
\newblock {\it Phys. Rev. D} {\bf 25} (1982), 1587-1604.

\bibitem{GR}
G. von ~Gehlen and R. Rittenberg.
\newblock $\Z_n$-symmetric quantum chains with infinite set of
conserved charges and $\Z_n$ zero modes.
\newblock {\it Nucl. Phys. B} {\bf 257} (1985), 351-370.

\bibitem{HandT}
B. ~Hartwig and P. ~Terwilliger.
\newblock The Tetrahedrom Algebra, the Onsager algebra, and the
$\mathfrak{sl}_2$ loop algebra.
\newblock {\it J. Algebra}, submitted.

\bibitem{iandt}
T. ~Ito and P. Terwilliger.
\newblock Tridiagonal pairs and the quantum affine algebra $U_q
(\widehat{sl}_2)$.
\newblock {\it Ramanujan J.}, to appear.

\bibitem{ITT}
T. ~Ito, K. ~Tanabe, and P. ~Terwilliger.
\newblock Some algebra related to P- and Q-polynomial association schemes.
\newblock {\it Codes and Association Schemes (Piscataway NJ, 1999)},
Amer. Math. Soc., Providence RI, 2000.

\bibitem{ito}
T. ~Ito and P. ~Terwilliger.
\newblock The shape of a tridiagonal pair.
\newblock {\it J. Pure Appl. Algebra,} {\bf 188} (2004) 145--160.

\bibitem{ons1}
L. ~Onsager.
\newblock Crystal statistics. I. A two-dimensional model with an
order-disorder transition.
\newblock {\it Phys. Rev. (2),} {\bf 65} (1944) 117--149.

\bibitem{P}
J.H.H. ~Perk.
\newblock Star-triangle relations, quantum Lax pairs, and higher
genus curves.
\newblock {\it Proceedings of Symposia in Pure Mathematics} {\bf
49} 341--354. Amer. Math. Soc., Providence, RI, 1989.

\bibitem{roan}
S. S. ~Roan.
\newblock Onsager's algebra, loop algebra and chiral Potts model.
\newblock Preprint MPI 91-70, Max Plank Institute for Mathematics,
Bonn, 1991.

\bibitem{ter2}
P. Terwilliger.
\newblock Two relations that generalize the $q$-Serre and Dolan-Grady
relations.
\newblock Proceedings of Nagoya 1999 Workshop on Physics and
Combinatorics, Nagoya, Japan 1999, World ScientificPublishing Co.,
Inc., River Edge, NJ, Providence, RI, 2000.

\bibitem{ter3}
P. ~Terwilliger.
\newblock Leonard pairs from 24 points of view.
\newblock {\it Rocky Mountain J. Math.} {\bf 32} (2)
(2002) 827--888.


\end{thebibliography}
\end{document}